\documentclass[11pt]{article}

\input{amssym.def}
\input{amssym}
\usepackage{eufrak}
\usepackage{color}

\newtheorem{theorem}{Theorem}

\newtheorem{remark}[theorem]{Remark}

\def\pa{\partial}
\def\h{\hbar}
\def\rh{r_\hbar}
\def\De{\Delta}

\def\AA{\cal A}
\def\BB{\cal B}

\def\Ah{{\cal A}_\hh}
\def\hh{\hbar}

\def\parh{\pa_{\rh}}

\def\de{\delta}

\def\ot{\otimes}

\def\C{{\Bbb C}}

\def\AA{{\cal A}}

\def\un{\underline}

\def\lhq{\ifmmode {\cal L}(q,\hbar)\else ${\cal L}(q,\hbar)$\fi}

\def\lqh{\ifmmode {\cal L}(q,\hbar)\else ${\cal L}(q,\hbar)$\fi}

\def\hpa{\hat{\pa}}
\def\hD{\hat{D}}
\def\DD{\mathrm{D}}

\def\hatt{\hat{\Theta}}

\def\la{{\lambda}}
\def\al{{\alpha}}

\def\be{\begin{equation}}
\def\ee{\end{equation}}

\begin{document}

\makeatletter
\renewcommand{\theequation}{{\thesection}.{\arabic{equation}}}
\@addtoreset{equation}{section} \makeatother

\title{Yang-Mills models on  Noncommutative background}

\author{\rule{0pt}{7mm} Dimitry
Gurevich\thanks{gurevich@ihes.fr}\\
{\small\it Institute for Information Transmission Problems}\\
{\small\it Bolshoy Karetny per. 19,  Moscow 127051, Russian Federation}\\
\rule{0pt}{7mm} Pavel Saponov\thanks{Pavel.Saponov@ihep.ru}\\
{\small\it
National Research University Higher School of Economics,}\\
{\small\it 20 Myasnitskaya Ulitsa, Moscow 101000, Russian Federation}\\
{\small \it and}\\
{\small \it
Institute for High Energy Physics, NRC "Kurchatov Institute"}\\
{\small \it Protvino 142281, Russian Federation}}

\maketitle

\begin{abstract}
In our previous publications we have developed some elements of Noncommutative calculus on the enveloping algebras of $A_m$ type, in particular, analogs of the partial
derivatives and de Rham complex were defined. Also, we introduced the notion of quantization with  Noncommutative configuration space and quantized  a few dynamical
models in this sense. In the current paper we propose a method of quantizing the Yang-Mills models in same  sense. 
\end{abstract}

{\bf AMS Mathematics Subject Classification, 2010:} 46L52, 81T13


\section{Introduction}

There are known numerous approaches  to constructing   dynamical models on noncommutative (NC) configuration spaces (more precisely, algebras).
The most known are models based on the algebras with a weak noncommutativity, defined via the following relations  between the generators
 $$
x_i\, x_j-x_i\, x_j=\theta_{ij},\quad \theta_{ij}\in \C.
$$
The product in the corresponding algebras is usually introduced via the famous Moyal formula.
In these algebras the partial derivatives $\pa_{x_i}$ can be defined in the classical way. It should be emphasized that defining analogs of the partial derivatives
in algebras with a more sophisticated noncommutativity  is not so easy. Nevertheless, we need them for considering models with nontrivial momenta.

In \cite{GPS2} we have introduced analogs of the partial derivatives on the algebras\footnote{The subscript $h$ means that we have introduced the multiplier $h$
in the standard $gl(m)$ Lie bracket. This allows us to treat the algebra $U(gl(m)_h)$ as a deformation (quantization) of the algebra $Sym(gl(m))$, corresponding
to the value $h=0$.} $U(gl(m)_h)$ and also on their super- and $q$-analogs. Below, these operators acting in the algebra $U(gl(m)_h)$ (or in its compact counterpart
$U(u(m)_h)$) will be called {\em quantum partial derivatives} (QPD). Note that the action of the QPD is subject to a Leibniz rule, which is a deformation of the usual one,
valid in the commutative algebra $Sym(gl(m))$.

By $q$-analogs of $U(gl(m)_h)$ we mean the so-called (modified) Reflection Equation algebras, associated with Hecke type braidings called Hecke symmetries.
In the present note we consider neither super- nor  $q$-analogs, the reader is referred to \cite{GS1}--\cite{GS3} for the corresponding constructions.

Making  use of the QPD we can define {\em quantum differential operators}, constructed according to the classical pattern
$$
\sum a_{i_1\dots\, i_k} \pa_{i_1}\dots \pa_{i_k}, \quad a_{i_1\dots\, i_k}\in U(u(m)_h),
$$
where $\pa_i$ stands for the QPD corresponding to $i$-th $gl(m)_h$ generator of the enveloping algebra $U(gl(m)_h)$. Moreover, the central step of our construction
consists in assigning such a quantum differential operator to its classical counterpart, corresponding to the value $h=0$ and acting in  $Sym(gl(m))$. This assignment
procedure, called {\em the quantization with noncommutative configuration space}, consists in replacing the usual derivatives by their quantum counterparts and the coefficients,
assumed to be elements of $Sym(gl(m))$, by their counterparts from $U(gl(m)_h)$\footnote{These counterparts arise from a map $\al_h: Sym(gl(m))\to U(gl(m)_h)$ which is
$GL(m)$-covariant and smoothly depends on $h$. Observe that such a map induces a new product on the algebra $Sym(gl(m))$ from the algebra  $U(gl(m)_h)$.}.

However, if the coefficients of a given operator are not polynomials,  we have to deal with an algebra larger than $U(gl(m)_h)$ and therefore we have to extend  the QPD
onto this larger algebra in order to  enlarge the class of possible solutions to the NC dynamical systems. In \cite{GS5} (also, see the reference therein) we dealt with the algebra
$U(u(2)_h)$ generated by the standard elements $x$, $y$, $z$ and $t$, playing the role of the spatial variables and time correspondingly, and introduced the notion of the quantum radius
$\rh=\sqrt{x^2+y^2+z^2+\h^2}$, $\h=h/2i$, which is an element of the central extension of this algebra. We succeeded in extending the QPD onto the quotient algebra
\be
\Ah=(U(su(2)_h)\ot \C[t, \rh])/\langle x^2+y^2+z^2+\h^2-\rh^2\rangle
\label{Ah}
\ee
and even onto some elements of the skew-field ${\BB}={\AA}[{\AA}^{-1}]$ with preserving the quantum Leibniz rule. Here, the notation $\langle J \rangle$ stands for the ideal generated
by a family $J$.

Note that in definition of the algebra $\Ah$ the polynomial factor $\C[t, \rh]$ can be replaced by the algebra of rational functions or Laurent series in two commuting generators $t$
and $\rh$. Below, we will use the notation $\C(t, \rh)$ for such a larger algebra. This algebra coincides with the center $Z(\Ah)$ of the algebra $\Ah$.

In the current paper we suggest a way of using the algebra $\Ah$ in quantization of the Yang-Mills models in the above sense. The problem is that  in our approach
the connection form (i.e. the Yang-Mills field)
of these models is represented by matrices with NC entries. As a result, the components of the Yang-Mills field strength lose the basic property of being skew-symmetric
with respect to  the permutation of the indexes $F_{\nu \mu}=-F_{\mu \nu}$. Below, we suggest a way of solving this problem by making use of the symmetrized form for the entries $F_{\mu\nu}$.

As an example we consider the 't Hooft-Polyakov monopole in the form of the Bogomol'niy equation and get a system describing an NC version of the so-called hedgehog. Since
the 't Hooft-Polyakov hedgehog is a spherically symmetric model, the central role in its description is played by the radial variable. In our NC version a similar role belongs to the quantum radius $\rh$.

The paper is organized as follows. In the next section we recall the construction of the quantum partial derivatives and exhibit the Leibniz rule for them. In Section 3 we introduce
the quantum radius and related operators. The NC version of the 't Hooft-Polyakov hedgehog is described in Section 4. In the last section we compare our version of NC geometry
with that from \cite{MS}  and explain why we consider our procedure of quantization with NC configuration spaces as a method of discretization.

\section{Quantum Leibniz rule }

It should be emphasized from the very beginning that it is very convenient  to write down all relations in the algebras we are dealing with in a matrix form. Thus, consider
the standard generators $l_i^j$, $1\leq i, j \leq m$ of the algebra $U(gl(m)_h)$ subject to the following commutation relations:
$$
l_i^j\, l_k^s-l_k^s\, l_i^j=h(l_i^s\,\de_k^j-l_k^j\, \de_i^s), \quad 1\leq i,j,k,l\leq m.
$$
Then, introducing the $m\times m$ matrix $L=\|l_i^j\|_{1\leq i, j \leq m}$, where the lower index numerates the rows and the upper one --- the columns, we can cast the above relations
in the matrix form:
$$
L_1\, L_2-L_2\, L_1=h (L_1\, P_{12}-P_{12}\, L_1).
$$
Hereafter, $P$ stands for the matrix of the usual flip and the lower index of a matrix indicates the way of its embedding into larger matrices. For example, for any $k\ge 1$ we
define
$$
L_k = I^{\otimes (k-1)}\otimes L\otimes I^{\otimes (N-k)},
$$
where $I$ is the $m\times m$ unit matrix, and a fixed number $N$ is assumed to be sufficiently large in order all our formulae make sense.

Now, consider a matrix $D$, whose entries $\pa_i^j$ play the role of the QPD corresponding to generators $l_j^i$ respectively, i.e. $\pa_i^j=\pa_{l_j^i}$. The QPD commute with each other, in the matrix
form this fact is expressed by the equality:
$$
D_1 D_2 =D_2D_1.
$$
The matrix form of relations between generators $\pa_k^l$ and $l_i^j$ reads:
\be
D_1 L_2=L_2 D_1+P_{12}+h D_1P_{12}.
\label{tri}
\ee
The map (denoted $\beta$), which sends the left hand side of (\ref{tri}) to its right hand side is called the permutation map, whereas the relations (\ref{tri}) themselves are called the permutation ones. Also, we assume that
$$\beta(\pa_i^j\ot 1)=1\ot \pa_i^j,\,\,\, \beta( 1\ot u)=u\ot 1,\,\, \forall \, u\in U(gl(m)_h).$$
The permutation map enables us to move any  polynomial  $p(\pa)$ through any polynomial $q(l)$  to the rightmost position. An important fact is that this procedure is in agreement
with the algebraic relations satisfied by the generators of the algebras involved.

Due to this fact, we can define an action of polynomials $p(\pa)$ onto these $q(l)$ as follows.  We apply the permutation map to the product $p(\pa)\ot q(l)$ moving all  QPD to
the rightmost position and then send to zero all terms containing at least one factor $\pa_k^l$. For details the reader is referred to \cite{GPS1, GPS2}. Note, that similarly to the
classical case ($h=0$) we have:
$$
\pa_i^j(1)=0\,\,\forall\, i,j,\quad \mathrm{and}\quad  \pa_i^j(l_k^l)=\de_i^l\, \de_k^j.
$$

Observe that the procedure of defining the QPD in $U(gl(m)_h)$ is a particular case of the  so-called quantum double $(A,B)$  of two associative algebras $A$ and
$B$ (see \cite{GS4}, where numerous types of quantum doubles are exhibited). In the case under consideration the algebra $A$ (respectively, $B$) consists of the polynomials in the QPD
$\pa_k^l$ (respectively in generators\footnote{The quantum double construction allows one to introduce  analogs of the partial derivatives on the Reflection Equation algebras as well.}
$l_i^j$).

Note that similarly to  the classical case the Leibniz rule for the QPD  can be expressed via a coproduct $\De:A\to A\ot A$. On the generators $\pa_i^j$ the map $\Delta$
acts as follows:
\be
\De(\pa_i^j)=\pa_i^j\ot 1_A+1_A\ot \pa_i^j +h\sum_k \pa_k^j\ot \pa^k_i:={\pa_i^j}_{\!(1)}\ot {\pa_i^j}_{\!(2)}.
\label{first}
\ee
Also, we put $\De(1_A)=1_A\ot 1_A$ and extend the coproduct on the whole algebra $A$ in the way converting it into a bialgebra.

\begin{remark}\rm
Observe that we did not succeed in finding an analogous coproduct for the QPD defined on the Reflection Equation algebras. On these algebras it is possible to define only a form of the Leibniz rule based on the permutation relations with subsequent canceling the terms containing the factors $\pa_i^j$. Whereas, a coproduct (\ref{first}) is very useful tool of the differential calculus in the algebras $U(gl(m)_h)$, since it enables us to find the action of any QPD on a product of elements from this algebra if we know its action on each factor.
\end{remark}

It is convenient to make the following shift of the generators of the algebra $A$:
\be
\hat{D}=D+\frac{I}{h} \quad \Leftrightarrow \quad \hpa_i^j=\pa_i^j+\frac{1}{h}\,\delta_i^j\,1_A.
\label{modi}
\ee
In terms of the matrix $\hD$ the permutation relations (\ref{tri}) acquires the following form:
\be
\hD_1 L_2=L_2 \hD_1+h \hD_1P_{12}.
\label{triii}
\ee
The action of the coproduct $\Delta$ (\ref{first}) on the generators $\hpa_i^j$ also becomes more simple:
\be
\De(\hpa_i^j)=h\sum_k \hpa_k^j\ot \hpa_i^k \quad \Leftrightarrow \quad \De(\hD^t)=h \hD^t\stackrel{.}{\ot} \hD^t,
\ee
where  $t$ denotes  the matrix transposition, and the notation $M\stackrel{.}{\ot} N$ stands for the matrix with entries $m_i^k\ot n_k^j$ for any two square matrices
$M=\|m_i^j\|$ and  $N=\|n_i^j\|$ of the same size.

Let us introduce the  matrix $\DD=h\hD^t$.  For this matrix the coproduct reads:
\be
\De(\DD)=\DD\stackrel{.}{\ot} \DD.
\label{Leiii}
\ee
Now we define the action $\hD(a)$ of the matrix $\hD$ on an arbitrary element $a\in U(gl(m)_h)$ in a  matrix form:
$$
\hD(a) = \|\hpa_i^j(a)\|, \quad \Rightarrow\quad \DD(a) = h\hD(a)^t.
$$
As a consequence of (\ref{Leiii}) we have:
\be
\DD(a b)=\DD(a) \DD(b),\quad \forall \,a, b \in U(gl(m)_h).
\label{L-rule-mat}
\ee

Thus, the matrix $\DD$ is a multiplicative (group-like) object with respect to the coproduct $\De$. Below we use the quantum Leibniz rule in the form (\ref{L-rule-mat}).
Besides, when extending the action of the QPD onto a larger algebra (for instance, on a central extension of $U(gl(m)_h)$ and on the corresponding skew-field), we recquire
the property (\ref{L-rule-mat}) of the matrix $\DD$ to be preserved. Thus, let $a\not=0$ be a element of the algebra $B=U(gl(m)_h)$, for which the values 
$\hpa_i^j(a)$ are known. Then, the  values $\hpa_i^j(a^{-1})$ can be found from the following relation:
$$
\DD(a^{-1})\DD(a)=\DD(a^{-1\,}a)=\DD(1_B)= I.
$$

In the case $m=2$ we are dealing with the matrices
\be
L =  \left(\!\!
\begin{array}{cc}
l_1^1&l_1^2\\
l_2^1&l_1^2
\end{array}
\!\!\right):=\left(\!\!
\begin{array}{cc}
a&b\\
c&d
\end{array}
\!\!\right),\qquad
\hD  = \left(\!\!
\begin{array}{cc}
\hat \partial_1^1&\partial^2_1\\
\partial_2^1&\hat\partial_2^2
\end{array}
\!\!\right)= \left(\!\!
\begin{array}{cc}
\hat \partial_a&\partial_c\\
\partial_b&\hat\partial_d
\end{array}
\!\!\right).
\label{gen-m}
\ee

Hereafter, we omit  the ``hat'' notation  for  off-diagonal elements of the matrix $\hD$ since they are  not changed under the shift (\ref{modi}):
$\hat \partial_i^j\equiv \partial_i^j$ if $i\not= j$.

Then the explicit form of the permutation relations (\ref{triii}) reads:
$$
\begin{array}{llllll}
\hpa_a a = a\,\hpa_a+h \,\hpa_a & \hspace*{2mm} &\hpa_a b = b\, \hpa_a + h\,\pa_c& \qquad \pa_c a = a\, \pa_c& \hspace*{2mm} & \pa_c b = b\, \pa_c\\
\rule{0pt}{4.5mm}
\hpa_a c = c\,\hpa_a & \hspace*{2mm} & \hpa_a d = d\, \hpa_a &\qquad \pa_c c = c\, \pa_c + h\, \hpa_a & \hspace*{2mm} & \pa_c d = d\, \pa_c + h\,\pa_c\\
\rule{0pt}{6mm}
\pa_b a = a\, \pa_b + h\,\pa_b & \hspace*{2mm} & \pa_b b = b\, \pa_b + h\,\hpa_d&\qquad \hpa_d a = a\, \hpa_d & \hspace*{2mm} & \hpa_d b = b\, \hpa_d\\
\rule{0pt}{4.5mm}
\pa_b c = c\,\pa_b & \hspace*{2mm} & \pa_b d = d\, \pa_b & \qquad \hpa_d c = c\, \hpa_d + h\,\pa_b & \hspace*{2mm} & \hpa_d d = d\, \hpa_d + h\,\hpa_d.
\end{array}
$$

To pass to the compact form $u(2)_\h $ of the algebra $gl(2)_\h$, we introduce the new set of generators
$$
t={{1}\over{2}}(a+d), \qquad x={{i}\over{2}}(b+c), \qquad y= {{1}\over{2}}(c-b), \qquad z={{i} \over{2}}(a-d)
$$
with the following  Lie brackets:
$$
[x, \, y]=h z,\qquad [y, \, z]=h x,\qquad[z, \, x]=h y,\qquad [t, \, x]=[t, \, y]=[t, \, z]=0.
$$

The QPD $\pa_t...\pa_z$ can be expressed via  these $\pa_a...\pa_d$ in the classical way. We leave exhibiting the permutation relations between the new generators and  the corresponding QPD
 to the reader. We only note that the shift (\ref{modi}) does not affect the QPD $\pa_x$, $\pa_y$, and $\pa_z$. So, only the derivative in $t$ becomes shifted:
$\hpa_t=\pa_t+\frac{2}{h}\, {\rm id}$. Its permutation relations with the generators of the algebra $u(2)_h$ are
$$
\hpa_t \,t - t\,\hpa_t = \frac{h}{2}\,\hpa_t,\quad \hpa_t \,x - x\,\hpa_t= - \frac{h}{2}\,\pa_x, \quad \hpa_t \,y - y\,\hpa_t=-\frac{h}{2}\,\pa_y, \quad \hpa_t \,z - z\,\hpa_t= - \frac{h}{2}\,\pa_z.
$$

Below, together with the  quantization parameter $h$ we will use its renormalized form $\h=h/2i$. Dealing with the algebra  $u(2)_h$, it is convenient to use the matrix:
\be
{\hatt}=
i\hh \left(\begin{array}{rrrr}
\hpa_t&\pa_x&\pa_y&\pa_z\\
-\pa_x&\hpa_t&-\pa_z&\pa_y\\
-\pa_y&\pa_z&\hpa_t&-\pa_x\\
-\pa_z&-\pa_y&\pa_x&\hpa_t
\end{array} \right)
\label{seven}
\ee
instead of $\DD$ defined above.

Note that the matrix  ${\hatt}$ is also multiplicative with respect to the coproduct $\De$:
$\De(\hatt)=\hatt\stackrel{.}{\ot}\hatt$.
Thus, the linear map
\be
\hatt: \, a\mapsto  \hatt(a)\quad  \forall\, a\in U(u(2)_h)
\label{mapp}
\ee
defines a map $U(u(2)_h)\to \mathrm{Mat}_4(U(u(2)_h))$, which preserves the product
\be
\hatt(ab)=\hatt(a) \hatt(b),\quad \forall\,a,b\in U(u(2)_h).
\label{mult}
\ee
Thus, we get a representation of the algebra $U(u(2)_h)$, in which any element of $U(u(2)_h)$ is represented by a $4\times 4$ matrix with entries from this algebra.

Now, as an example, we compute the action $\pa_z(x  y)$. First, we note that  the QPD $\pa_t$, $\pa_x$, $\pa_y$ and $\pa_z$ act on the generators $t$, $x$, $y$ and $z$
in the classical way\footnote{It is also valid for  the unit $1\in U(u(2)_h)$.}, and then apply the coproduct
$$
\De(\pa_z)=i\h\,(\hpa_t\ot \pa_z+\pa_z\ot  \hpa_t+\pa_x\ot  \pa_y- \pa_y\ot  \pa_x),
$$
which follows from the mentioned property $\Delta\hatt = \hatt\stackrel{.}{\ot}\hatt$. As a result we get:
$$
\pa_z(x y)=i\h(\hpa_t(x)\,\pa_z(y)+\pa_z(x)\,\hpa_t(y)+\pa_x(x)\,\pa_y(y)- \pa_y(x)\,\pa_x(y))=i\h=\frac{h}{2}.
$$
Thus, $\pa_z(x y)=0$ only in the limit $\h=0$.

\section{Quantum radius and quantum radial derivative}

As can be verified by a straightforward calculation, the matrix $L$ generating the algebra $U(u(2)_h)$
$$
L = \left(\!\!\!
\begin{array}{rcr}
t-iz&\,&-ix-y\\
\rule{0pt}{5mm}
-ix+y&&t+iz
\end{array}\!\!
\right)
$$
satisfies a matrix polynomial identity of the form $p(L)=0$,
where
\be
p(s)=s^2-(2t+\h)\,s+(t^2+x^2+y^2+z^2+2i\hh\,t)\, 1.
\label{poly}
\ee
Note, that all coefficients of the polynomial $p(s)$ belong to the center $Z(U(u(2)_h)$. Analogous identities are valid  for the generating matrices of the algebras
$U(gl(m)_h)$, their super-analogues and the Reflection Equation algebras. All of them are called Cayley-Hamilton identities. 

Let $\mu_i$, $i=1,2$ be the roots of the polynomial (\ref{poly}), that is $p(s) = (s-\mu_1)(s-\mu_2)$, or equivalently:
$$
\mu_1+\mu_2 = 2t+\h,\qquad \mu_1\mu_2 = t^2+x^2+y^2+z^2+2i\hh\,t.
$$
 We consider the central extension $U(u(2)_h)[\mu_1,\mu_2]$ of the initial algebra $U(gl(2)_h)$ and call the quantities $\mu_1$ and $\mu_2$ {\em the eigenvalues} of the matrix $L$.

In \cite{GS5} we extended the action of the QPD on $\mu_i$, $i=1,2$, with preserving the quantum Leibniz rule (\ref{L-rule-mat}). Here, we only present the results of the QPD action
on the quantum radius $\rh = \frac{\mu}{2 i}$, $\mu=\mu_1-\mu_2$, which belongs to the center of the extended algebra $U(u(2)_h)[\mu_1,\mu_2]$. Applying the map (\ref{mapp}) to
the quantum radius we get (see \cite{GS5} for detail):
$$
\hatt(\rh)=\frac{\rh^2+\h^2}{\rh}\, I+\frac{i\, \h}{\rh}\,
\left(\!\!
\begin{array}{cccc}
 0 & x&y&z\\
 -x&0&-z&y\\
-y&z&  0 & -x\\
-z&y& x & 0
\end{array}
\!\!\right).
$$
The matrix elements of $\hatt(\rh)$ define the actions of all QPD onto the quantum radius. Namely, we have:
\be
\pa_t \rh=-\frac{i\hh}{\rh},\qquad \pa_x\rh=\frac{x}{\rh},\qquad \pa_y\rh=\frac{y}{\rh},\qquad \pa_z\rh=\frac{z}{\rh}.
\label{pervy}
\ee
It should be emphasized that three last formulae in (\ref{pervy}) look like the classical ones. As for the first one, it turns into the classical result at the limit $\h=0$.
Note also, that all considerations do not depend on the enumeration order of the eigenvalues $\mu_i$, which is arbitrary.

In what follows we need a generalization of the QPD action on any rational function (or a formal series) $f( \tau, \rh)$, where $\tau=-i t$ is the renormalized
time\footnote{The writing $f(\tau, \rh)$ is not ambiguous since the generators $\tau$ and $\rh$ commute with each other.}.
The corresponding formulae were found in \cite{GS2}:
$$
\pa_\tau(f(\tau, \rh))=\frac{f(\tau+\h,\rh+\h)(\tau+\h)+f(\tau+\h,\rh-\h)(\tau-\h)-2\rh f(\tau, \rh)}{2\rh\h},
$$
$$
\pa_x(f(\tau, \rh))=\frac{x}{\rh}\frac{f(\tau+\h,\rh+\h)-f(\tau+\h,\rh-\h)}{2\h}.
$$
The formulae for  $\pa_y(f(\tau, \rh))$ and $ \pa_z(f(\tau, \rh))$ are similar. It is worth pointing out that using the 
generator $\tau$ allows us to avoid the complex numbers in the above formulae.

Now, we define {\em the quantum radial derivative} $\parh$ as an operator, acting on functions $f(\tau,\rh)$ by the following rule:
\be
\parh (f(\tau, \rh))=\frac{f(\tau+\h, \rh+\h)-f(\tau+\h, \rh-\h)}{2\h}, \qquad \forall\, f(\tau, \rh)\in \C(\tau, \rh).
\label{rhh}
\ee
It should be emphasized that the quantum radial derivative affects the time generator. However, this property disappears at the classical limit $\h=0$. Besides, if $\h\not=0$,
$\parh$ is a difference operator with the respect to the both generators $\rh$ and $\tau$.

By using the quantum radial derivative we can represent the above formulae  in the spirit of  the classical relations:
$$\pa_x f(\tau, \rh)=\frac{x}{\rh}\, \parh f(\tau, \rh),\quad \pa_y f(\tau, \rh)=\frac{y}{\rh} \,\parh f(\tau, \rh),\quad \pa_z f(\tau, \rh)=\frac{z}{\rh} \,\parh f(\tau, \rh).
$$
However, so far the operator $\parh$ is defined only on the functions $f(\tau, \rh)$.

In order to extend  the action of the derivative $\parh$ onto other elements of the algebra $\Ah$ we  use the so-called  (quantum) Laplacian  $\Delta: U(u(2)_h)\to U(u(2)_h)$:
$$
\De=\pa^2_x+\pa^2_y+\pa^2_z.
$$

The action of the operators $\De$ on polynomials in $\tau$ and  $\rh$ turns out to be (see \cite{GS2} for detail):
\begin{eqnarray*}
\De(f(\tau, \rh))= \frac{f(\tau+\h, \rh+2\h)+f(\tau+\h, \rh-2\h)-2 f(\tau+\h, \rh)}{4\h^2}&&\\
&&\hspace*{-35mm}+\,\frac{f(\tau+\h, \rh+2\h)-f(\tau+\h, \rh-2\h)}{2\h\, \rh}.
\end{eqnarray*}
Taking into account definition (\ref{rhh}) we have:
$$
\De(f(\tau, \rh))=\frac{1}{\rh} \, \parh^2(\rh f(\tau, \rh)).
$$
This relation  is similar to the classical one:
$$
\De(f(t, r))=\frac{1}{r}\,\pa^2_r (r  f(t, r))=\pa^2_r  f(t, r)+\frac{2}{r}\,\pa_r (f(t, r)).
$$

Besides, we consider the operator $Q=x\, \pa_x+y\, \pa_y+\, z\,\pa_z$. In the classical setting (at $\h=0)$ this operator equals $Q\,\rule[-1.5mm]{0.4pt}{3mm}_{\,\h=0}  = r\,\pa_r$.
In the current NC setting this relation becomes (see \cite{GS2}):
$$
Q(f(\tau, \rh))=\frac{\rh^2-\h^2}{\rh}\, \parh (f(\tau, \rh)).
$$
On requiring this connection between the operators $Q$ and $\parh$ to be valid on the arbitrary element $a\in \Ah$ by definition, we can extend the action of the derivative
$\parh$ onto the whole algebra $\Ah$, defined by   (\ref{Ah}).

As  examples we compute the actions $\pa_x (x f(\tau, \rh))$ and $\pa_x (y f(\tau, \rh))$. Thus, we have
$$
\pa_x (x\,f(\tau, \rh))=i\h(\hpa_t(x)\cdot \pa_x(f)+\pa_x(x)\cdot \hpa_t(f)+\pa_y(x)\cdot \pa_z(f)-\pa_z(x)\cdot \pa_y(f))=\frac{x^2}{\rh}\parh(f)+f+\h\pa_\tau(f),
$$
$$
\pa_x (y\,f(\tau, \rh))=i\h(\hpa_t(y)\cdot \pa_x(f)+\pa_x(y)\cdot \hpa_t(f)+\pa_y(y)\cdot \pa_z(f)-\pa_z(y)\cdot \pa_y(f)=\frac{yx+i\h z}{\rh}\parh(f).
$$
These formulae will be  essentially used in the next section.

\section{Yang-Mills models and their quantization}

Given a simple Lie algebra with a basis $T_i$ and structure constants $c_{ij}^k$, i.e. $[T_i, T_j]=c_{ij}^k\, T_k$, the corresponding Yang-Mills field is by definition
a connection form $A_\mu^{\,\, i}$ such that the corresponding covariant derivative is
$$
\nabla_\mu=I\, \pa_\mu+A_\mu^{\,\, i}\, T_i.
$$

The  field strength $F_{\mu \nu}$ is defined via the commutator of the operators $\nabla_\mu$ and $\nabla_\nu$:
$$
F_{\mu\nu} = [\nabla_\mu,\nabla_\nu] = F_{\mu \nu}^{\,\,\, i}T_i,
$$
where
\be
F_{\mu \nu}^{\,\,\, i}=\pa_\mu \, A_\nu^{\,\, i}-\pa_\nu \, A_\mu^{\,\, i}+ A_\mu^{\,\, k}\, A_\nu^{\,\, l}\, c_{kl}^{\,\,i}.
\label{stre}
\ee
Since the structure constants $c_{kl}^{\,\,i}$ are skew-symmetric in the lower indexes, we have
$$
F_{\nu \mu}^{\,\,\, i}=-F_{\mu \nu}^{\,\,\, i}.
$$
It should be emphasized that this property is valid since the quantities $A_\mu^{\,\, i}$ are functions and consequently commute with each other. However, this property fails if we assume
$A_\mu^{\,\, i}$ to be elements of an NC algebra.

In order to restore the skew-symmetry of $F_{\mu\nu}$ when the components $A_\mu^{\,\, i}$ do not commute with each other, we have to replace the product
$A_\mu^{\,\, k} A_\nu^{\,\, l}$ by its symmetric part:
$$
\frac{1}{2}(A_\mu^{\,\, k} A_\nu^{\,\, l}+A_\nu^{\,\, l} A_\mu^{\,\, k}).
$$
Below, for this symmetric part we use the notation $\un{A_\mu^{\,\, k} A_\nu^{\,\, l}}$. Thus, instead of (\ref{stre}) we have the following formula defining components of the field
strength in an NC case:
$$
F_{\mu \nu}^{\,\,\, i}=\pa_\mu  A_\nu^{\,\, i}-\pa_\nu A_\mu^{\,\, i}+ \un{A_\mu^{\,\, k} A_\nu^{\,\, l}}\, c_{kl}^{\,\,i}.
$$

Besides, we realize the same operation for any term containing a product of elements from the NC algebra in question. For instance, in the 't Hooft-Polyakov monopole model
the products $A_\mu^{\,\,k}\varphi^i$, where $\{\varphi^i\}_{ i=1,2,3}$ is the triplet of scalar fields, will be replaced by
$$
\un{A_\mu^{\,\,k}\,\varphi^i}=\frac{1}{2}(A_\mu^{\,\,k}\varphi^i+\varphi^i A_\mu^{\,\,k}).
$$

Now, let us recall the 't Hooft-Polyakov monopole model under  the Bogomol'nyi equation form (see \cite{K})
\be
F_{\mu \nu}^{\,\,\, i}=\varepsilon_{\mu \nu}^{\,\,\,\la}\nabla_\la\, \varphi^i.
\label{dyn-eq}
\ee

Here, the Lie algebra involved is  $su(2)$. By fixing in it the usual orthonormal basis formed by the Pauli matrices (up to  normalizing factors), we have that the structure constants
$c_{kl}^{\,\,i}=\varepsilon_{kl}^{\,\,i}$ are skew-symmetric with respect to all indexes (we assume that  $\varepsilon_{12}^{\,\,3}=1$). The corresponding metric, which serves for raising
and lowering indexes, is Euclidean $\de_{ij}$. So, we do not distinguish the low and upper indexes.

Next, we assume the elements $A_\mu^{\,\,k}$ and $\varphi^i$ to be elements of the enveloping algebra $U(su(2))$ and, taking in the mind the above modifications, we
are looking for solutions of the Bogoml'nyi equation by using the following ansatz:
\be
A_\mu^{\,\, i}=\varepsilon_\mu^{\,\,\, ij} x_j\, W(\rh),\quad \varphi^i=x^i\, F(\rh).
\label{ans}
\ee
Here, we use the natural identification  $x_1 = x$, $x_2 = y$, $x_3 = z$. The notations for the QPD are modified correspondingly.

To get the equations on structure functions $W$ and $F$ it is sufficient to use two equations from the system (\ref{dyn-eq}):
$$
F_{12}^{\,\,\, 1}= \pa_x A_2^1-\pa_y A_1^1+\un{A_1^2 A_2^3} - \un{A_1^3 A_2^3}=\nabla_3 \varphi^1=\pa_z \varphi^1+\un{A_3^2\, \varphi^3}-\un{A_3^3\varphi^2}
$$
and
$$
F_{12}^{\,\,\, 3}= \pa_x A_2^3-\pa_y A_1^3+\un{A_1^1 A_2^2}-\un{A_1^2 A_2^1}=\nabla_3 \varphi^3=\pa_z \varphi^3+\un{A_3^1\, \varphi^2}-\un{A_3^2\varphi^1}.
$$
Due to a specific form of ansatz (\ref{ans}) the other equations lead to the same results.

On substituting expressions (\ref{ans}) into the above dynamical equations we get the system of equations on the functions $W$ and $F$:
$$
-\pa_x(z W) +\un {zx}\, W^2=\pa_z(xF)-\un{xz} \,FW,
$$
$$
\pa_x(xW)+\pa_y(yW)+z^2\, W^2=\pa_z(zF)+y^2 FW+x^2 \,FW.
$$
After symmetrizing the underlined terms and cancelling the common   factor $zx+xz$, the first equation reduces to the form:
\be
-\frac{1}{\rh}\,{\parh W}+W^2=\frac{1}{\rh}\,{\parh F}-FW.
\label{odi}
\ee
The second equation can be transformed to the following relation:
$$
\frac{x^2+y^2}{\rh}\, \parh W+2W +2\h \,\partial_\tau W + z^2\,W^2=\frac{z^2}{\rh}\,\parh F+F+\h\, \partial_\tau F +(x^2+y^2) FW.
$$
Upon  replacing  the sum $x^2+y^2$ by $\rh^2-\h^2-z^2$ and taking into account equation (\ref{odi}) we get the final form of the second equation:
\be
\frac{\rh^2-\h^2}{\rh} \,\parh W+ 2W+2\h \,\partial_\tau W=F+\h\, \partial_\tau F+(r^2-\h^2) FW.
\label{dva}
\ee

Thus, we get two equations for two unknown functions $W(\rh)$ and $F(\rh)$. Let us point out the main difference between
these equations and their classical limit: besides the
presence of new terms with $\h$-factors, the QPD $\parh$ is actually a {\it difference} operator. It would be interesting to find a nontrivial solution of the system (\ref{odi})--(\ref{dva}).

\section{Concluding remarks}

In this note we do not consider the NC analogs of the differential forms. In other words, we restrict ourselves to deformation of the classical Weyl-Heisenberg (WH) algebra, generated by commutative algebra $Sym(gl(m))$ and the corresponding partial derivatives. A different  deformation of the WH algebras with the same configuration space (and some other similar ones)
 has been  introduced in \cite{MS}. However, the permutation relations of the corresponding deformed WH algebra in that paper are defined via series in the initial derivatives. This leads to more complicated formulae, which do not allow to introduce analogs of quantum radius and related objects.

Let us mention two other advantages of our approach. First, it is based on a quite simple Leibniz rule which fairly simplifies all calculations, especially for spherically symmetric models.
We plan to quantize with NC configuration space some other physically significant models in our future publications.  

Second, from our approach it becomes clear that the usual partial derivatives become difference operators  on the center of the algebra $\Ah$. In this sense we can say that a passage
to our NC configuration space and the corresponding QPD is a tool of discretization. In the classical approach one deals with commutative algebras of functions. Whereas we deal with
the enveloping algebra, which has completely different representation theory: its finite-dimensional irreducible representations are enumerated by the spin and the dynamic can be
considered as a passage from one representation to another one. 
If the algebra $U(u(2)_h)$ is represented in the space of spin $n$, the quantum radius takes the value $\rh=(2n+1)\h$
provided the generators $x$, $y$ and $z$ are represented by Hermitian operators. Also, we assume the parameter $\h$  to be real and positive.

For a more detailed discussion on the discretization of models in the indicated sense the reader is referred to \cite{GS2}. In particular, in that paper we noticed that some operators
(for instance, the radial part of the Laplace operator) can be restricted on a lattice, defined in the space generated by $\rh$ and $\tau$.

It is also interesting to study the action of the QPD onto the components of the algebra $\Ah$, which are complementary to the center. 
Observe  that the algebra $U(su(2)_h)$ can be decomposed in the direct sum of subspaces:
$$
U(su(2)_h)=\bigoplus_{k=0}^\infty \left( Z(U(su(2)_h))\ot V^{(k)}\right),\,\, k=0,1,2...
$$
where $V^{(k)}\subset U(su(2)_h)$ is the spin $k$ component.  For $k=0$ we have  $V^{(0)}=\C$. The highest weight element in the complexification of the component $V^{(k)}$
is $b^k$, where $b=l_1^2$ (see (\ref{gen-m})). The center $Z(U(su(2)_h))$ is generated by the Casimir element $x^2+y^2+z^2$. This entails that the algebra $\Ah$ is spanned by
the elements $f(\tau,\rh) \ot u$, $u\in V^{(k)}$, $\forall\,k\ge 0$.

As was shown in \cite{GS2}, the results of applying  the QPD $\pa_t,\dots,\pa_z$ to the components $V^{(k)}$ are classical. This fact can be checked directly by showing that the
actions of the derivatives $\pa_a,\, \pa_b,\, \pa_c, \pa_d$ onto the elements $b^k$, $\forall\, k\ge 0$, are classical. In this sense our partial derivatives become quantum indeed only on the elements $f(\tau, \rh)\in  Z(\Ah)$.

As a final observation we mention the paper \cite{IS}, where our QPD have been used for a direct transferring to the enveloping algebras the shift argument method, developed
in the Poisson setting for constructing commutative subalgebras. The problem of quantizing such subalgebras has been formulated by E.Vinberg. A method of solving this problem was
proposed in \cite{T}. However, the method developed in \cite{IS} and making use of our QPD is more simple and straightforward.


\begin{thebibliography}{PPP}

\bibitem[GPS1]{GPS1} Gurevich D., Pytov P., Saponov P., Braided Differential operators on Quantum algebras,  J. of Geometry and
Physics 61 (2011), pp. 1485--1501.

\bibitem[GPS2]{GPS2} Gurevich D., Pytov P., Saponov P.,    Braided Weyl algebras and differential calculus on $U(u(2))$, J. of Geometry and
Physics 62 (2012), pp. 1175--1188.

\bibitem[GS1]{GS1} Gurevich D., Saponov P., Braided algebras and their applications to Noncommutative Geometry,  Advances in Applied Mathematics  51, (2013) 228--253.

\bibitem[GS2]{GS2} Gurevich D., Saponov P., Noncommutative Geometry and dynamical models on $U(u(2))$ background, Journal of Generalized Lie theory and Applications, 9:1 (2015).

\bibitem[GS3]{GS3} Gurevich D., Saponov P., Quantum geometry and  quantization on $U(u(2))$ background. Noncommutative Dirac monopole, J. of Geometry and Physics, 106 (2016), pp.   87-97.

\bibitem[GS4]{GS4} Gurevich D., Saponov P., Doubles of Associative algebras and their Applications, Phys. of Particles and Nuclei Letters 17 N5 (2020) 774--778.

\bibitem[GS5]{GS5} Gurevich D., Saponov P.,  Noncommutative geometry on central extension of $U(u(2))$,  J. of  Mathematical Physics 64 (2023).

\bibitem[IS]{IS}  Ikeda Y.,  Sharygin G.I., The argument shift method in universal enveloping algebra $U(gl_d)$,  arXiv:2307.15952.

\bibitem[K]{K} Katanaev M. Spherically symmetric 't Hooft-Polyakov monopole, Eur. Phys. J. C (2021) 81:825.

\bibitem[MS]{MS}  Meljanac S., \v{S}koda Z.,  Leibniz rules  for enveloping algebras in symmetric ordering, arXiv:0711.0149.

\bibitem[T]{T} Tarasov A.A., On some commutative subalgebras in the universal enveloping algebra of
the Lie algebra $gl(n, C)$, Sb. Math. 191 (2000), 1375–1382.

\end{thebibliography}
\end{document}